\documentstyle[12pt,leqno,amssymb,graphics]{amsart}

\newtheorem{thm}{Theorem}[section]
\newtheorem{lemma}[thm]{Lemma}
\newtheorem{prop}[thm]{Proposition}

\theoremstyle{definition}
\newtheorem{dfn}[thm]{Definition}
\theoremstyle{remark}

\newtheorem{example}[thm]{Example}

\begin{document}

\newcommand{\ct}{\cite}
\newcommand{\pr}{\protect\ref}
\newcommand{\su}{\subseteq}
\newcommand{\pa}{{\partial}}
\newcommand{\im}{{Imm(F,\E)}}
\newcommand{\hf}{{1 \over 2}}
\newcommand{\Q}{{\mathbb Q}}
\newcommand{\R}{{\mathbb R}}
\newcommand{\Z}{{\mathbb Z}}
\newcommand{\G}{{\mathbb G}}
\newcommand{\F}{{\mathbb F}}

\newcommand{\h}{\widehat}

\newcommand{\X}{{\mathbb X}}
\newcommand{\Y}{{\mathbb Y}}
\newcommand{\U}{{\mathbb U}}
\newcommand{\E}{{{\math R}^3}}
\newcommand{\I}{{\mathrm{Id}}}
\newcommand{\4}{{\mathcal{H}}}
\newcommand{\C}{{\mathcal{C}}}
\newcommand{\cc}{{\C_0}}
\newcommand{\1}{{(1)}}

\newcommand{\ce}{{\mathcal{C}^{\text{ev}}}}
\newcommand{\co}{{\mathcal{C}^{\text{od}}}}

\newcommand{\cce}{{\mathcal{C}^{\text{ev}}_0}}
\newcommand{\cco}{{\mathcal{C}^{\text{od}}_0}}

\newcommand{\xe}{{\X^{\text{ev}}}}
\newcommand{\xo}{{\X^{\text{od}}}}
\newcommand{\ye}{{\Y^{\text{ev}}}}
\newcommand{\yo}{{\Y^{\text{od}}}}
\newcommand{\ev}{{\text{ev}}}
\newcommand{\od}{{\text{od}}}

\newcommand{\we}{{\W^{\text{ev}}}}
\newcommand{\wo}{{\W^{\text{od}}}}

\newcommand{\hc}{{H_1(F,\C)}}
\newcommand{\ak}{{ \{ a_k \}  }}   
\newcommand{\bk}{{ \{ b_k \}  }}   
\newcommand{\tb}{{ \Leftrightarrow }} 
\newcommand{\bn}{{ \leftrightarrow }} 

\newcounter{numb}

\title{Order one invariants of spherical curves}
\author{Tahl Nowik}
\address{Department of Mathematics, Bar-Ilan University, 
Ramat-Gan 52900, Israel}
\email{tahl@@math.biu.ac.il}
\date{October 10, 2007}
\urladdr{http://www.math.biu.ac.il/$\sim$tahl}

\begin{abstract}
We give a complete description of all order 1 invariants of spherical curves. 
We also identify the subspaces of all $J$-invariants and $S$-invariants, 
and present two equalities satisfied by any spherical curve.
\end{abstract}

\maketitle

\section{Introduction}\label{intro}

A spherical or planar curve is an immersion of $S^1$ in $S^2$ or $\R^2$ respectively.
The study of invariants of planar and spherical curves has been initiated by Arnold in \ct{a1}--\ct{a3},
where he presented the three basic order 1 invariants $J^+,J^-,St$.
Various explicit formulas for these and other invariants appear in \ct{cd},\ct{p1},\ct{p3},\ct{s},\ct{vi}.
The study of order 1 and higher order invariants has in general split into the study of $J$-invariants,
as in \ct{g1},\ct{g2},\ct{l},\ct{p2}, 
and $S$-invariants as in \ct{k},\ct{t},\ct{va}, 
where $J$-invariants are invariants which are unchanged when passing a triple
point, and $S$-invariants are invariants which are unchanged when passing a tangency.
But as we shall see for the case of spherical curves, the space of invariants spanned by the order 1
$J$- and $S$-invariants is much smaller than the full space of order 1 invariants.

In this work we give a complete presentation of all order 1 invariants of spherical curves.
Since the vector space in which the invariants take their values is arbitrary, the presentation
is in terms of a universal order 1 invariant.
We also identify the subspaces of $J$-invariants and $S$-invariants, 
and present the formulas for Arnold's three invariants which stem from our presentation.
In the final section we present two peculiar equalities satisfied by any spherical curve.

Analogous study for immersions of surfaces in 3-space has been carried out in \ct{n1}--\ct{n5}.
In fact, our present universal order 1 invariant of spherical curves involves similar ingredients to 
those involved in the invariant for immersions of surfaces in 3-space introduced in \ct{n4}.

The structure of the paper is as follows. In Section \pr{state} we 
present the basic definitions, construct our invariant $F$ of spherical curves, and state our results.  
In Section \pr{abs} we construct an abstract invariant, as opposed to the explicit invariant $F$, 
and show that the abstract invariant is universal (Theorem \pr{absf}).
In Section \pr{prf} we prove our main result, namely, that the explicit invariant $F$ is universal
(Theorem \pr{main}).
In Section \pr{siji} we construct universal invariants for the spaces of $S$-invariants 
and $J$-invariants. We also present formulas for Arnold's invariants $J^+,J^-,St$ for spherical curves.
In the concluding Section \pr{oee} we discuss the way in which all our work divides between the two regular
homotopy classes of spherical curves, and present two equalities satisfied by any spherical curve (Theorem \pr{fin}).

\section{Definitions and statement of results}\label{state}
By a \emph{curve} we will always mean an immersion $c:S^1 \to S^2$.   
Let $\C$ denote the space of all curves. The space $\C$ has two connected components, that is, two
regular homotopy classes, which we denote $\ce,\co$. 
A curve $c$ is in $\ce$ (respectively $\co$)
if when deleting a point $p$  from $S^2$ not in the image of $c$, 
the curve $c$, as a curve in $S^2 - \{ p \} \cong \R^2$ has even winding number (respectively odd winding number). 
A curve will be called \emph{stable} if its only self intersections are transverse double points.
A stable curve $c$ is in $\ce$ (respectively $\co$) 
iff the number of double points of $c$ is odd (respectively even). 

The generic singularities a curve may have are either a tangency of first order between two strands, 
which will be called a $J$-type singularity, or three strands meeting at a point, 
each two of which are transverse, which will be called an $S$-type singularity.
Singularities of type $J$ and $S$ appear in Figures \pr{c2} and \pr{c3}. 
A generic singularity can be resolved in two ways, 
and there is a standard way for considering one resolution positive, and the other negative,
as defined in \ct{a1}.  
We denote by $\C_n \su \C$ ($n\geq 0$) the space of all curves which have precisely 
$n$ generic singularity points (the self intersection being elsewhere stable).
In particular, $\C_0$ is the space of all stable curves. 
An invariant of curves is a function $f:\C_0 \to W$,
which is constant on the connected components of $\C_0$, and where $W$ in this work 
will always be a vector space over $\Q$.

Given a curve $c\in \C_n$, with singularities located at $p_1,\dots,p_n \in S^2$, 
and given a subset $A\su \{p_1,\dots,p_n\}$,
we define $c_A \in C_0$ to be the stable curve obtained from $c$ by resolving all singularities
of $c$ at points of $A$ into the 
negative side, and all singularities not in $A$ into the positive side.
Given an invariant $f:\C_0 \to W$ we define the ``$n$th derivative'' of $f$ to be the function
$f^{(n)}:\C_n \to W$ defined by
$$f^{(n)}(c)=\sum_{ A \su \{p_1,\dots,p_n\} } (-1)^{|A|} f(c_A)$$
where $|A|$ is the number of elements in $A$.
An invariant $f:\C_0\to W$ is called \emph{of order $n$} if 
$f^{(n+1)}(c)=0$ for all $c\in \C_{n+1}$.
The space of all $W$ valued invariants on $\C_0$ of order $n$ is denoted $V_n=V_n(W)$.
Clearly $V_n \su V_m$ for $n \leq m$.
In this work we give a full description of $V_1$. We will construct a ``universal'' 
order 1 invariant, by which we mean the following:

\begin{dfn}\label{uni}
An order 1 invariant $\h{f} :\C_0 \to \h{W}$ will be called \emph{universal}, if
for any $W$ and any order 1 invariant $f:\C_0 \to W$, there exists a unique linear map 
$\phi: \h{W} \to W$ such that $f = \phi \circ \h{f}$. In other words, for any $W$, the natural map
$Hom_\Q(\h{W},W) \to V_1(W)$ given by $\phi \mapsto \phi \circ \h{f}$ is an isomorphism.
\end{dfn}

\begin{dfn}\label{gen2}
Let $D$ be a 2-disc. An immersion $e:[0,1] \to D$ will be called a \emph{simple arc} if
\begin{enumerate}
\item $e(0),e(1) \in \pa D$.
\item $e$ is transverse to $\pa D$.
\item The self intersections of $e$ are transverse double points.
\end{enumerate}
\end{dfn}

\begin{dfn}\label{di}
Let $D$ be an oriented 2-disc, and let $e$ be a simple arc in $D$.
\begin{enumerate}
\item For a double point $v$ of $e$
we define $i(v) \in \{1,-1\}$, where $i(v)=1$ if the orientation at $v$ given by the two tangents to $e$ at $v$,
in the order they are visited, coincides with the orientation of $D$. Otherwise $i(v)=-1$.
\item We define the index of $e$, $i(e)\in\Z$ by $i(e) = \sum_v i(v)$ 
where the sum is over all double points $v$ of $e$.
\end{enumerate}
\end{dfn}

The following is easy to see, and is the reason for defining $i$:

\begin{lemma}\label{li}
If $e,e'$ are two simple arcs in $D$ with the same initial and final points and the same initial and
final tangents, then $i(e)=i(e')$ iff $e$ and $e'$ are regularly homotopic in $D$, keeping the initial and final
points and tangents fixed. 
\end{lemma}

Define $\X$ to be the vector space over $\Q$ with basis all symbols $X_{a,b}$ where $(a,b) \in \Z^2$
is an ordered pair of integers.
Let $\Y$ be the vector space over $\Q$ with basis all symbols $Y_d$ where $d \in \Z$.
We construct two invariants $f^X : \C_0 \to \X$ and $f^Y: \C_0 \to \Y$ in the following way.

For $c \in \C_0$ let $v$ be a double point of $c$, and let $u_1, u_2$ be the two tangents at  
$v$ ordered by the orientation of $S^2$.
Let $U$ be a small neighborhood of $v$ and let $D=S^2 - U$. 
Now $c|_{c^{-1}(D)}$ defines two simple arcs $c_1,c_2$ in $D$, ordered so that
the tangent $u_i$ leads to $c_i$, $i=1,2$. They will be called the exterior arcs of $c$.
We denote $a(v) = i(c_1)$ and $b(v) = i(c_2)$ 
where the orientation on $D$ is that restricted from $S^2$.
We define $f^X : \C_0 \to \X$ as follows:
$$f^X(c) = \sum_v X_{a(v),b(v)}$$ 
where the sum is over all double points $v$ of $c$.

If $c \in \C_0$ then the image of $c$ divides $S^2$ into complementary regions (which are open discs). 
To each such region $R$ we attach an integer $d(R)$ as follows. 
Take $p \in R$ and view $c$ as a curve in $S^2 - \{ p \} \cong \R^2$.
Let $d(R)$ be \emph{minus}
the winding number of $c$ as a curve in $S^2 - \{ p \}$ where the 
orientation on $S^2 - \{ p \}$ is that restricted from $S^2$.
We define $f^Y: \C_0 \to \Y$ as follows: $$f^Y(c) = \sum_R Y_{d(R)}$$
where the sum is over all complementary regions $R$ of $c$.

\begin{dfn}\label{F}
Let $F:\C_0 \to \X \oplus \Y$ be the invariant given by $F(c) = f^X(c) + f^Y(c)$.
\end{dfn}

\begin{figure}
\scalebox{0.8}{\includegraphics{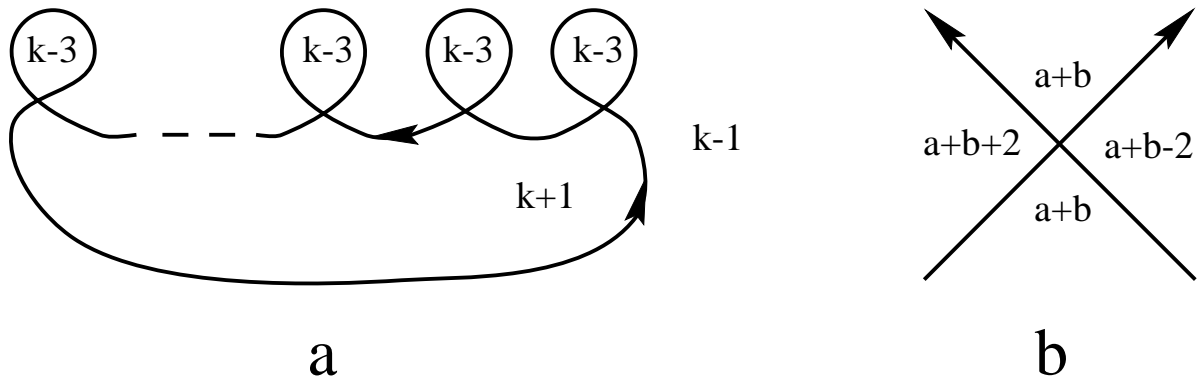}}
\caption{a: The curve $\Gamma_k$. b: $d(R)$ in the vicinity of a double point of type $(a,b)$.}\label{c1}
\end{figure}

\begin{example}
For each $k \geq 0$ let $\Gamma_k$ be the curve with $k$ double points appearing in Figure \pr{c1}a. 
In particular $\Gamma_0$ is the embedded circle and $\Gamma_1$ is the figure eight curve. 
We have $\Gamma_k \in \co$ for $k$ even, and $\Gamma_k \in \ce$ for $k$ odd. 
The integers marked in the figure are $d(R)$ of the various complementary regions.
For each of the $k$ double points we have $a(v)=0$ and $b(v)=k-1$, and so
we obtain $F(\Gamma_k) = k X_{0,k-1} + k Y_{k-3} + Y_{k-1} + Y_{k+1}$.
\end{example}

It is not hard to show that the pair of indices $a(v),b(v)$ of a double point $v$, determine $d(R)$
for the four regions adjacent to $v$, as stated in Figure \pr{c1}b (recall that $d(R)$ is \emph{minus} the 
corresponding winding number).

We define the following six linear maps $\psi_i:\X\oplus\Y \to \Q$, $i=1,\dots,6$ by stating their 
value on the basis $\{ X_{a,b},Y_d \}$, and where $\psi_i=0$ on each basis element which is not explicitly mentioned in 
its definition. 
\begin{enumerate}
\item $\psi_1(Y_d)=1$ for $d$ odd. $\psi_1(X_{a,b})=-1$ for $a+b$ odd. 
\item $\psi_2(Y_d)=1$ for $d$ even. $\psi_2(X_{a,b})=-1$ for $a+b$ even. 
\item $\psi_3(Y_d)=d$ for $d$ odd. $\psi_3(X_{a,b})=-(a+b)$ for $a+b$ odd. 
\item $\psi_4(Y_d)=d$ for $d$ even. $\psi_4(X_{a,b})=-(a+b)$ for $a+b$ even. 
\item $\psi_5(Y_d) = \frac{1}{d}$ for $d$ odd. 
$\psi_5(X_{a,b}) = \frac{4(a-b+1) - (a+b)^2}{(a+b)\left( (a+b)^2 - 4 \right)}$ for $(a+b)$ odd.
\item $\psi_6(Y_0)=1$. $\psi_6(X_{a,b}) = b-a-1$ for $a+b=0$. $\psi_6(X_{a,b}) = \frac{a-b}{2}$ for $a+b=\pm 2$.

\end{enumerate}

\begin{dfn}\label{eq} 
Let $W$ be a vector space over $\Q$, and let $\varphi_1,\dots,\varphi_n:W \to \Q$ be linear maps. 
We denote by $W^{\varphi_1,\dots,\varphi_n}$ the subspace of $W$ determined by the $n$ 
equations $\varphi_1=\varphi_2=\cdots=\varphi_n=0$.
\end{dfn}

Our main result, Theorem \pr{main}, will be proved in Section \pr{prf}, stating: 
$F(\C_0) \su (\X \oplus \Y)^{\psi_3,\psi_4,\psi_5,\psi_6}$, and 
$F:\C_0 \to (\X \oplus \Y)^{\psi_3,\psi_4,\psi_5,\psi_6}$ is a universal order 1 invariant.

In Section \pr{siji} we will construct universal invariants for the class of order 1
$S$-invariants, $J$-invariants, and for their combinations which we name $SJ$-invariants. 
In Section \pr{oee} we will explain how all our constructions split between $\ce$ and $\co$, and will 
also note the two equalities expressing the fact that $\psi_5,\psi_6$ vanish on the image of $F$.

\section{An abstract universal order one invariant}\label{abs}

In this section we will construct an ``abstract'' universal order 1 invariant, as opposed to the ``concrete''
invariant described in the previous section. 
This will be an intermediate step in proving that the concrete
invariant is universal. In addition, the existence and properties of the concrete invariant will aid us 
in proving properties of the abstract invariant.

Two curves $c,c' \in \C_1$ will be called \emph{equivalent} 
if there is an ambient isotopy of $S^2$ bringing a neighborhood $U$ of 
the singular point of $c$ onto a neighborhood $U'$ of the singular point of $c'$, such that the configuration
near the singular point precisely matches and such that the exterior arcs in $D = S^2-U$ are 
regularly homotopic in $D$. By Lemma \pr{li} a pair of exterior arcs are regularly homotopic
in $D$ iff their indices are the same. So we have:

\begin{prop}\label{eq1}
Two curves $c,c' \in \C_1$ are equivalent
iff  they have the same singularity configuration and the same corresponding indices of exterior arcs.
\end{prop}

The following is clear from the definition of order 1 invariant:

\begin{lemma}\label{eqv}
Let $f : \C_0 \to W$ be an invariant, then $f$ is of order 1 iff for any two equivalent $c,c' \in \C_1$,
$f^{(1)}(c) = f^{(1)}(c')$.
\end{lemma}

\begin{figure}
\scalebox{0.8}{\includegraphics{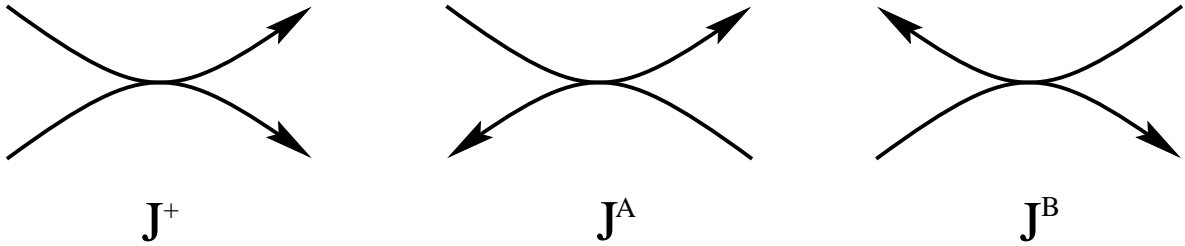}}
\caption{Singularities of type $J^+$, $J^A$, $J^B$.}\label{c2}
\end{figure}

We will attach a symbol to each equivalence class of curves in $\C_1$ as follows. For $J$ type singularities there are 
three distinct configurations which we name $J^+,J^A,J^B$, see Figure \pr{c2}.
The $J^A,J^B$ singularities are symmetric with respect to a $\pi$ rotation, which interchanges the two exterior arcs,
and so, by Proposition \pr{eq1},
the equivalence class of a $J^A$ or $J^B$ singularity is characterized by a symbol $J^A_{a,b}$ and $J^B_{a,b}$
($a,b \in \Z$), where $a,b$ is an \emph{unordered} pair, registering the indices of the two exterior arcs.
The $J^+$ configuration, on the other hand, is not symmetric, and so characterized by a symbol $J^+_{a,b}$ with 
$a,b$ an \emph{ordered} pair, with say, $a$ corresponding to the lower strand in Figure \pr{c2}.
(We will shortly see that this ordering may however be disregarded.) 

\begin{figure}
\scalebox{0.8}{\includegraphics{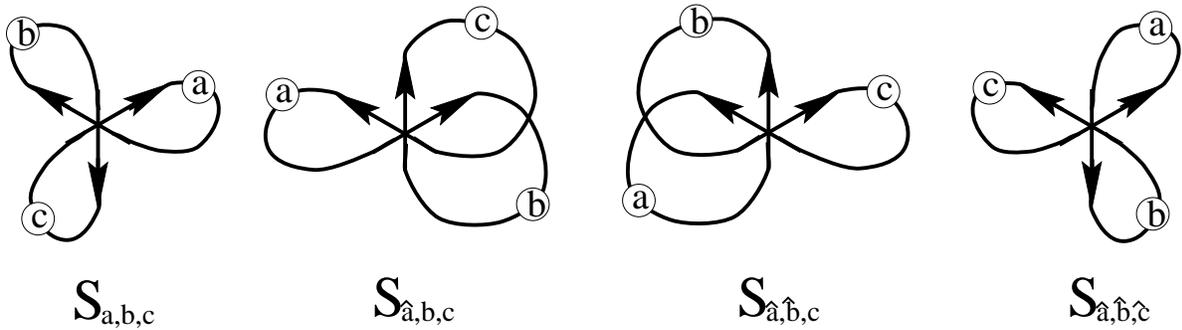}}
\caption{Singularities of type $S_{a,b,c}$, $S_{\h{a},b,c}$, $S_{\h{a},\h{b},c}$, $S_{\h{a},\h{b},\h{c}}$.}\label{c3}
\end{figure}

As to $S$ type singularities, there are four types, as seen in Figure \pr{c3}.
The distinction between them will be incorporated into the
way we register the indices of the exterior arcs. We will have a cyclicly ordered triple of integers 
$a,b,c$ registering the indices of the three exterior arcs, in the cyclic
order they appear along $S^1$, and each may appear with or without a hat, according to the following
rule. For given exterior arc $e$ with index $a$, let $u_1$ be the initial tangent of $e$ and $u_2$ the initial tangent
of the following segment. Then if $u_2$ is pointing to the right of $u_1$ then $a$ will appear with a hat,
and if $u_2$ is pointing to the left of $u_1$ then $a$ will appear unhatted. Since the ordering is cyclic, this gives
four types of $S$ symbols, $S_{a,b,c} , S_{\h{a},b,c} , S_{\h{a},\h{b},c}, S_{\h{a},\h{b},\h{c}}$, which
correspond to the four types of $S$ singularities.

Let $\F$ denote the vector space over $\Q$ with basis all the above symbols, $J^+_{a,b}$, $J^A_{a,b}$, $J^B_{a,b}$, 
$S_{a,b,c}$, $S_{\h{a},b,c}$, $S_{\h{a},\h{b},c}$, $S_{\h{a},\h{b},\h{c}}$.
If $\gamma$ is a generic path in $\C$, that is, a path passing only generic singularities,
then we denote by $v(\gamma) \in \F$ the sum of symbols of the singularities it passes, each added with $+$ or $-$
sign according to whether we pass it from its negative side to its positive side, or from its positive side to its 
negative side, respectively. Let $N \su \F$ be the subspace generated by all elements $v(\gamma)$ obtained from 
all possible generic \emph{loops} $\gamma$ in $\C$ (i.e. \emph{closed} paths), and let $\G = \F / N$. 

For an order 1 invariant $f:\C_0 \to W$, since $f^{(1)}$ coincides on equivalent curves in $\C_1$, it induces a
well defined linear map $f^{(1)}:\F \to W$. 

The following is clear:
\begin{lemma}\label{vgm}
If $\gamma$ is a generic path in $\C$, from $c_1$ to $c_2$, then $f^\1 ( v(\gamma) ) = f(c_2) - f (c_1)$. 
\end{lemma}

From Lemma \pr{vgm} it follows that $f^\1$ vanishes on the generators of $N$,
so it also induces a well defined linear map $f^{(1)} : \G \to W$.

Let $\h{\G} = \G \oplus \Q a_0 \oplus \Q a_1$,
where $a_0,a_1$ are two new vectors. We define an order 1 invariant $\h{f} : \C_0 \to \h{\G}$ as follows:
Denote $\C_n^\ev = \C_n \cap \ce$ and $\C_n^\od = \C_n \cap \co$, and recall $\Gamma_1 \in \C_0^\ev$,
$\Gamma_0 \in \C_0^\od$, are the figure eight curve, and embedded circle, respectively. 
For any $c \in \C_0^\ev$ there is a generic path $\gamma$ 
from $\Gamma_1$ to $c$. Let $\h{f}(c)=a_1 + v(\gamma) \in \h{\G}$. 
Similarly for $c \in \C_0^\od$ there is a generic path $\gamma$ 
from $\Gamma_0$ to $c$. Let $\h{f}(c)=a_0 + v(\gamma) \in \h{\G}$. 
By definition of $N$, 
$v(\gamma) \in \G$ is indeed independent of the choice of path $\gamma$.
From Lemma \pr{eqv}
it is clear that $\h{f}$
is an order 1 invariant, and we will now see that it is \emph{universal}: 

\begin{thm}\label{absf}
$\h{f} : \C_0 \to \h{\G}$ is a universal order 1 invariant. 
\end{thm}

\begin{pf}
For an order 1 invariant $f: \C_0 \to W$, define the linear map $\phi_f:\h{\G}\to W$
by $\phi_f|_{\G} = f^{(1)}$ and $\phi_f(a_i)=f(\Gamma_i)$, $i=0,1$.
We claim $\phi_f \circ \h{f} = f$ and that $\phi_f$ is the unique linear map satisfying this property.
Indeed, let $c \in \C_0^\ev$ and let $\gamma$ be a generic path from $\Gamma_1$ to $c$.
Then by Lemma \pr{vgm} 
$f(c) = f(\Gamma_1) + f^\1( v(\gamma) ) = \phi_f(a_1) + \phi_f ( v(\gamma)) = \phi_f ( \h{f} (c))$. 
Similarly this is shown for $c \in \C_0^\od$.

For uniqueness, it is enough to show that $\h{f} ( \C_0 )$ spans $\h{\G}$.
We have $a_i = \h{f}(\Gamma_i)$ for  $i=0,1$ so it remains to show that
$\G \su span \h{f} ( \C_0 )$.
Indeed for any generating symbol $T$ of $\G$, $T$ is the 
the difference $\h{f}(c)-\h{f}(c')$ for two
curves $c,c' \in \C_0$, namely, the two resolutions of a curve in $\C_1$ 
whose symbol is $T$.
\end{pf}

The following is also clear:

\begin{prop}\label{iso}
Let $f: \C_0 \to W$ be an order 1 invariant.
If $\phi_f: \h{\G} \to W$ (appearing in the proof of Theorem \pr{absf})
is an isomorphism, then $f:\C_0 \to W$
is also a universal order 1 invariant. 
\end{prop}

We now present six specific subfamilies of the set of generators of $N$.
We will eventually see (concluding paragraph of Section \pr{prf}) that these elements 
in fact span $N$.

\begin{figure}
\scalebox{0.8}{\includegraphics{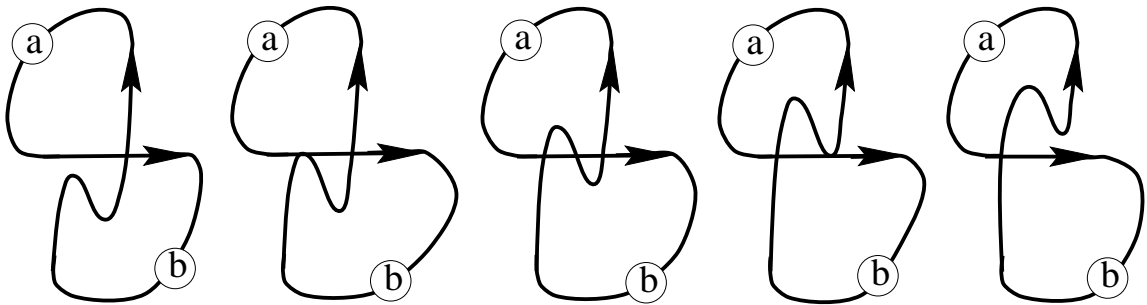}}
\caption{$J^+_{a,b} - J^+_{b,a} = 0$}\label{c4}
\end{figure}

\begin{figure}
\scalebox{0.8}{\includegraphics{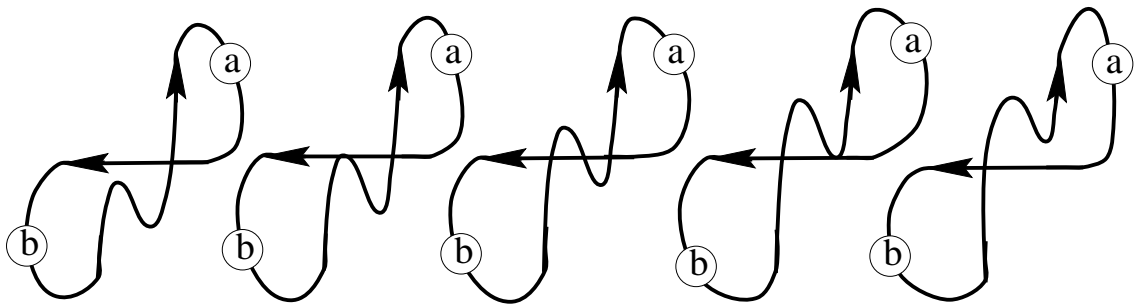}}
\caption{$J^B_{a+1,b} - J^A_{a,b-1} = 0$}\label{c5}
\end{figure}

The first two families of loops appear in Figures \pr{c4},\pr{c5}. 
In each white disc in a figure, there is assumed to be a simple arc of index given by the label on the disc.
The element obtained from Figure \pr{c4} is $J^+_{a,b} - J^+_{b,a}$, 
that is, in $\G$ we have the relation $J^+_{a,b} - J^+_{b,a} = 0$,
and so from now on we can simply regard the indices of $J^+_{a,b}$ as being \emph{un-ordered} 
(as is true by definition for the $J^A$ and $J^B$ symbols, 
and as opposed to the indices of $X_{a,b}$ which are ordered).
The relation obtained from Figure \pr{c5} is $J^B_{a+1,b} - J^A_{a,b-1} = 0$, or after index shift:
$J^B_{a,b} = J^A_{a-1,b-1}$.

\begin{figure}
\scalebox{0.8}{\includegraphics{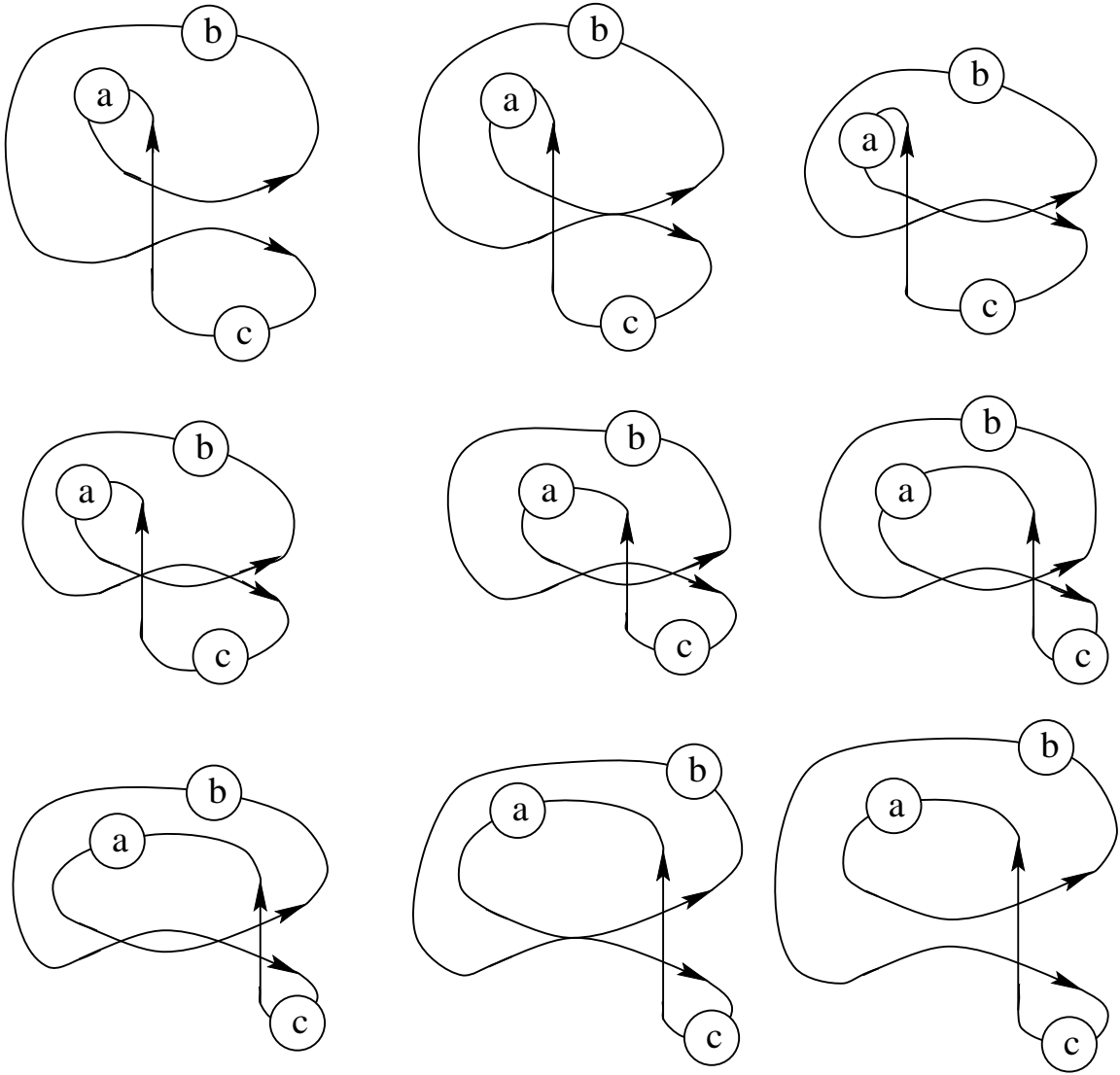}}
\caption{$J^+_{c+a-1,b} - S_{\h{a},b,c} + S_{\h{a},\h{b},c} - J^+_{c+a+1,b} =  0$}\label{c6}
\end{figure}

\begin{figure}
\scalebox{0.8}{\includegraphics{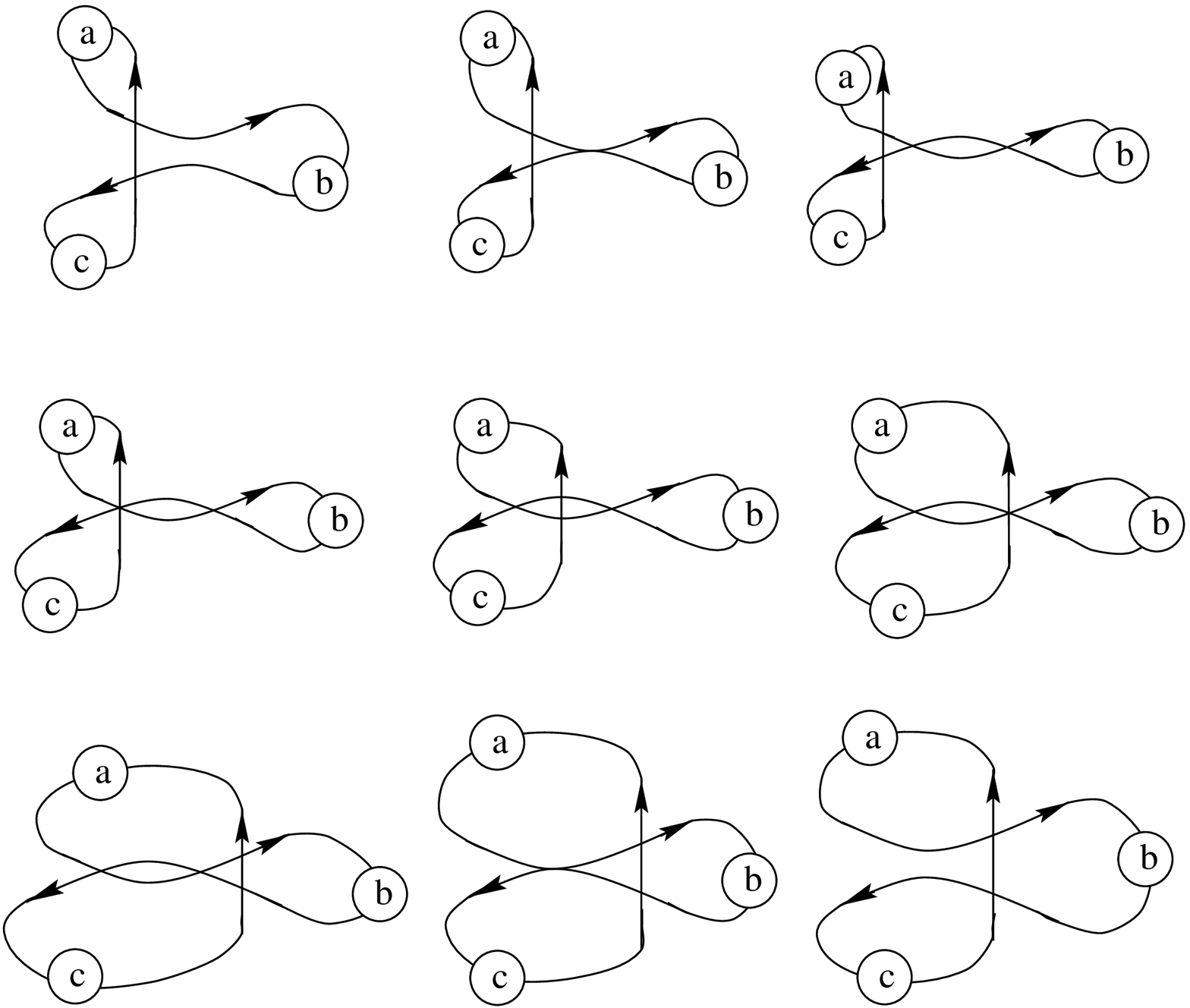}}
\caption{$J^A_{c+a-2,b} + S_{\h{a},\h{b+1},\h{c}} - S_{\h{a},b,\h{c}} - J^A_{c+a,b} =  0$}\label{c7}
\end{figure}

Our next two families of loops appear in Figures \pr{c6},\pr{c7}. 
The relations obtained are:
$J^+_{c+a-1,b} - S_{\h{a},b,c} + S_{\h{a},\h{b},c} - J^+_{c+a+1,b} =  0$ and 
$J^A_{c+a-2,b} + S_{\h{a},\h{b+1},\h{c}} - S_{\h{a},b,\h{c}} - J^A_{c+a,b} =  0$.
     
When taking the mirror image of such figure, 
the effect is that the indices are negated, $J^A$ and $J^B$ are interchanged,
and the ``hat status'' of each index is reversed.
So, the mirror image of the last two relations gives (after replacing $-a,-b,-c$ by $a,b,c$):
$J^+_{c+a+1,b} - S_{a,\h{b},\h{c}} + S_{a,b,\h{c}} - J^+_{c+a-1,b} =  0$ and 
$J^B_{c+a+2,b} + S_{a,b-1,c} - S_{a,\h{b},c} - J^B_{c+a,b} =  0$.
After replacing each $J^B_{a,b}$ with $J^A_{a-1,b-1}$ in the last relation, 
and an index shift, we obtain:
$J^A_{c+a+1,b} + S_{a,b,c} - S_{a,\h{b+1},c} - J^A_{c+a-1,b} =  0$.
We write these last four relations 
(coming from Figures \pr{c6},\pr{c7} and their mirror images)
in the following more convenient form:

\begin{prop}\label{relJS} 
The following relations hold in $\G$:
\begin{enumerate}
\item   $S_{\h{a},b,c} -  S_{\h{a},\h{b},c}              = J^+_{c+a-1,b} -  J^+_{c+a+1,b}$
\item   $S_{a,b,\h{c}} -  S_{a,\h{b},\h{c}}              = J^+_{c+a-1,b} -   J^+_{c+a+1,b}$
\item   $S_{\h{a},b,\h{c}} -  S_{\h{a},\h{b+1},\h{c}}    = J^A_{c+a-2,b} -    J^A_{c+a,b}$ 
\item   $S_{a,b,c} -  S_{a,\h{b+1},c}                    = J^A_{c+a-1,b} -    J^A_{c+a+1,b}$
\end{enumerate}
\end{prop}

We will have several occasions to use the following observation:

\begin{lemma}\label{sk}
The following two sets of relations are equivalent:
\begin{enumerate}
\item 
\begin{enumerate}
\item $S_{\h{a},b,c} =   S_{\h{a},\h{b},c} $
\item $S_{a,b,\h{c}} =   S_{a,\h{b},\h{c}} $
\item $S_{\h{a},b,\h{c}} =   S_{\h{a},\h{b+1},\h{c}} $ 
\item $S_{a,b,c} =   S_{a,\h{b+1},c} $
\end{enumerate}
\item Denoting $k = a+b+c$:
\begin{enumerate}
\item $S_{a,b,c} = S_{\h{k+1},0,0}$
\item $S_{\h{a},b,c} = S_{\h{a},\h{b},c} = S_{\h{k},0,0}$
\item $S_{\h{a},\h{b},\h{c}} = S_{\h{k-1},0,0}$
\end{enumerate}
\end{enumerate}

\end{lemma}

\begin{pf}
Relations 2 clearly imply relations 1. For the converse, 
alternatingly using 1a,1b we have 
$S_{\h{a},b,c} = S_{\h{a},\h{b},c} = S_{a,\h{b},c} = S_{a,\h{b},\h{c}} = S_{a,b,\h{c}}=S_{\h{c},\h{a},b}$,
that is, for given cyclic triple $a,b,c$, one or two hats at any location all give equal symbols.
From this and from 1d we have $S_{\h{a},b,c} = S_{a,\h{b},c} = S_{a,b-1,c} = S_{\h{a+1},b-1,c}$.
From all the above we obtain 2b. Now from 2b and 1c,1d, we obtain 2a,2c.
\end{pf}

\begin{dfn}We define the following two subsets of the generating set of $\F$:
\begin{enumerate}
\item $A^J = \{ J^+_{a,b} \}_{a,b\in\Z} \cup \{ J^A_{a,b} \}_{a,b\in\Z}$
\item $A^S = \{S_{\h{a},0,0}\}_{a \in\Z}$
\end{enumerate}
\end{dfn}

We will eventually see (Theorem \pr{c}(1)) that $A^J \cup A^S$ is a basis for $\G$. At this stage we show:

\begin{lemma}\label{gen}
$A^J \cup A^S$ is a spanning set for $\G$.
\end{lemma}

\begin{pf}
From the relation $J^B_{a,b}=J^A_{a-1,b-1}$ it follows that $A^J$ spans 
all symbols of type $J$. It follows from Proposition \pr{relJS} that in $\G / span A^J$ the relations 
of Lemma \pr{sk} hold, and so $A^S$ spans $\G / span A^J$.
\end{pf}

\section{Proof of main result}\label{prf}

Returning to our invariant $F:\C_1 \to \X \oplus \Y$ we now compute $F^{(1)}$.
It is easy to see that indeed the value of $F^\1$ depends only on the symbol of a
given curve $c \in \C_1$, which by Lemma \pr{eqv} proves  
that indeed $F$ is an order 1 invariant. For the sake of completeness, 
we find the value of $F^\1$ on all symbols, though we will need it only for   
the generators $A^J \cup A^S$ of $\G$.

\begin{figure}
\scalebox{0.8}{\includegraphics{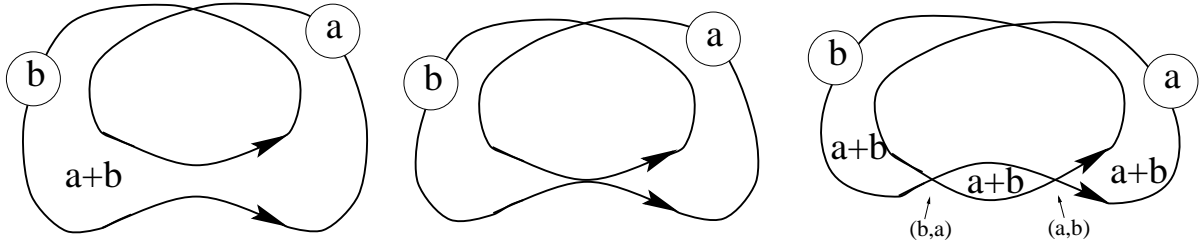}}
\caption{$F^{(1)}(J^+_{a,b})=X_{a,b} + X_{b,a} + 2Y_{a+b}$}\label{c8}
\end{figure}

\begin{figure}
\scalebox{0.8}{\includegraphics{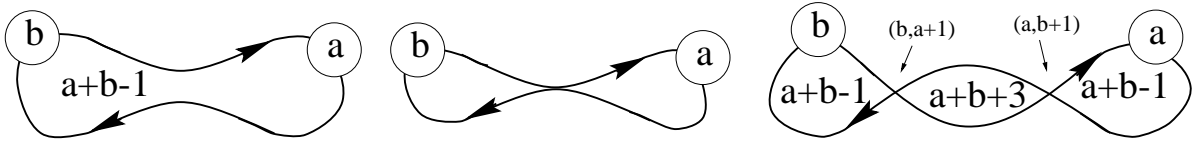}}
\caption{$F^{(1)}(J^A_{a,b})=X_{a,b+1} + X_{b,a+1} + Y_{a+b-1} + Y_{a+b+3}$}\label{c9}
\end{figure}

\begin{figure}
\scalebox{0.8}{\includegraphics{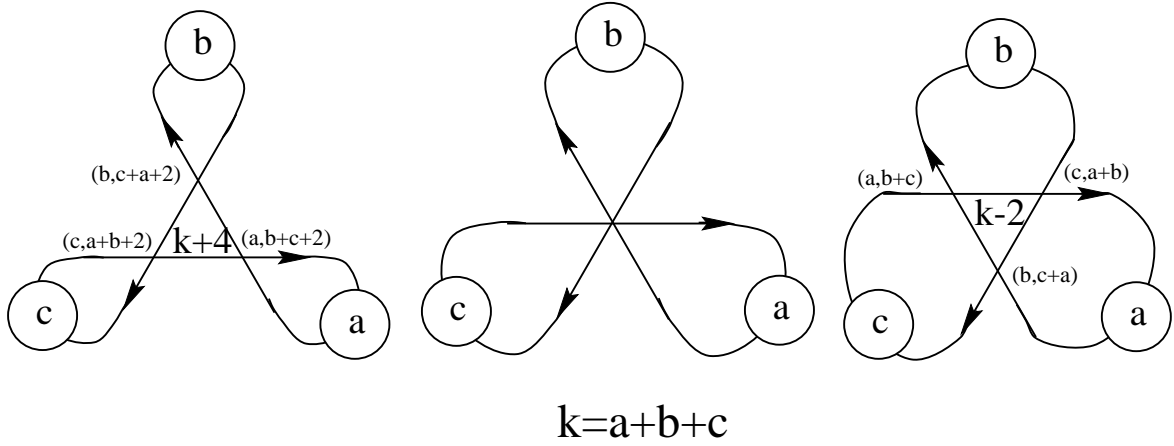}}
\caption{$F^\1(S_{a,b,c})=   -X_{a,b+c+2}-X_{b,c+a+2}-X_{c,a+b+2} +X_{a,b+c} +X_{b,c+a} +X_{c,a+b} 
  -Y_{a+b+c+4}+Y_{a+b+c-2}$}\label{c10}
\end{figure}

By Figures \pr{c8},\pr{c9} we have 
\begin{itemize}
\item $F^{(1)}(J^+_{a,b})=X_{a,b} + X_{b,a} + 2Y_{a+b}$ 
\item $F^{(1)}(J^A_{a,b})=X_{a,b+1} + X_{b,a+1} + Y_{a+b-1} + Y_{a+b+3}$
\end{itemize}
It follows from this and the relation $J^B_{a,b} = J^A_{a-1,b-1}$ that 
\begin{itemize}
\item $F^\1(J^B_{a,b})=X_{a-1,b}+X_{b-1,a}+Y_{a+b-3} + Y_{a+b+1}$
\end{itemize}
By Figure \pr{c10} we have:
\begin{itemize}
\item $F^\1(S_{a,b,c})=   -X_{a,b+c+2}-X_{b,c+a+2}-X_{c,a+b+2} +X_{a,b+c} +X_{b,c+a} +X_{c,a+b}    \\
  -Y_{a+b+c+4}+Y_{a+b+c-2}$
\end{itemize}
It follows from the above values of $F^\1$ and the relations of Proposition \pr{relJS} that  
\begin{itemize}
\item $F^\1(S_{\h{a},b,c})=   -X_{c,a+b+1}-X_{b+c-1,a}-X_{b,c+a+1}+X_{b+c+1,a}+X_{b,c+a-1}+X_{c,a+b-1} \\
  -Y_{a+b+c+1}+Y_{a+b+c-1}$
\item $F^\1(S_{\h{a},\h{b},c})=   -X_{c,a+b+1}-X_{c+a-1,b}-X_{b+c-1,a}+X_{b+c+1,a}+X_{c,a+b-1}+X_{c+a+1,b} \\
  -Y_{a+b+c-1}+Y_{a+b+c+1}$
\item $F^\1(S_{\h{a},\h{b},\h{c}})=   -X_{b+c-2,a}-X_{c+a-2,b}-X_{a+b-2,c}+X_{a+b,c}+X_{b+c,a}+X_{c+a,b} \\
  -Y_{a+b+c-4}+Y_{a+b+c+2}$
\end{itemize}

We now perform the following change of basis for $\X \oplus \Y$.
Define 
$$\Phi_{a,b}=\frac{a-b}{2}Y_{a+b-2} + (b-a-1)Y_{a+b} + \frac{a-b}{2}Y_{a+b+2}$$ 
and let $Z_{a,b}= X_{a,b} - \Phi_{a,b}$. 
Then we claim $\{Z_{a,b},Y_d\}$ is a basis for $\X \oplus \Y$. 
Indeed the linear map from $\X \oplus \Y$ to itself taking 
$X_{a,b} \mapsto Z_{a,b}$, $Y_d \mapsto Y_d$ which is given
in the $\{ X_{a,b} , Y_d \}$ basis by $X_{a,b} \mapsto X_{a,b} - \Phi_{a,b}$, $Y_d \mapsto Y_d$,
has inverse given by $X_{a,b} \mapsto X_{a,b} + \Phi_{a,b}$, $Y_d \mapsto Y_d$. 
(Since $\Phi_{a,b} \in \Y$ it is fixed by these maps.)

Substituting 
$X_{a,b}=Z_{a,b}+\Phi_{a,b} $
in the above expressions for $F^\1$, on the generators $A^J \cup A^S$, 
we obtain:
\begin{itemize}
\item $F^{(1)}(J^+_{a,b})=Z_{a,b} + Z_{b,a}$. 
\item $F^{(1)}(J^A_{a,b})=Z_{a,b+1} + Z_{b,a+1}$.
\item $F^\1(S_{\h{a},0,0})= - 2Z_{0,a+1}-Z_{-1,a}+Z_{1,a}+ 2Z_{0,a-1} \\
 -\frac{1}{2}(a-3)Y_{a-3} +\frac{3}{2}(a-1)Y_{a-1} - \frac{3}{2}(a+1)Y_{a+1} + \frac{1}{2}(a+3)Y_{a+3}$.
\end{itemize}

We see that the change of basis has simplified the formulas for $F^{(1)}(J^+_{a,b}), F^{(1)}(J^A_{a,b})$.
It also simplifies the formulas for $\psi_1,\dots,\psi_6$:

\begin{prop}\label{vanish}
The linear maps $\psi_1,\dots,\psi_6 : \X \oplus \Y \to \Q$ of Section \pr{state} 
all vanish on each $Z_{a,b}$.
\end{prop}

\begin{pf}
This technical verification is best done by defining
new linear maps $\psi'_i : \X \oplus \Y \to \Q$ $i=1,\dots,6$, 
via the $\{ Z_{a,b} , Y_d \}$ basis, by $\psi'_i(Z_{a,b}) = 0$ and $\psi'_i(Y_d) = \psi_i(Y_d)$. 
Then substitute 
$X_{a,b}=Z_{a,b}+ \Phi_{a,b} $
into the formulas for $\psi'_i$ and check that the formulas for $\psi_i$ are obtained.
\end{pf}

The following lemma contains the linear algebra portion of our analysis.

\begin{lemma}\label{inj}
The images of $A^J \cup A^S$ under  $F^\1:  \G \to \X \oplus \Y$ are independent, and span the subspace
$(\X \oplus \Y)^{\psi_1, \dots , \psi_6}$.
\end{lemma}

\begin{pf}
We perform all computations in the $\{Z_{a,b},Y_d\}$ basis.
We first show that $F^\1(A^J)$ are independent in $\X \oplus \Y$, and span the subspace 
$Z = span \{ Z_{a,b} \}$.
That is, we will show that the set $\{ Z_{a,b} + Z_{b,a} : a \geq b\} \cup \{ Z_{a,b+1} + Z_{b,a+1} : a \geq b\}$ 
is a basis  for $Z$.
For this it is enough to show that for any fixed $n \in \Z$, 
$\{ Z_{a,b} + Z_{b,a} : a \geq b, a+b =n \} \cup \{ Z_{a,b+1} + Z_{b,a+1} : a \geq b, a+b = n-1 \}$ 
is a basis for $span \{ Z_{a,b} : a+b =n \}$. For this it is enough to show that 
for every $k  \leq \frac{n}{2}$, the $n-2k+1$ dimensional subspace 
$span\{ Z_{a,b} :  a \geq k,  b \geq k , a+b =n\}$ is spanned by the $n-2k+1$ elements
$\{ Z_{a,b} + Z_{b,a} :  a \geq b \geq k , a+b =n\} \cup  \{ Z_{a,b+1} + Z_{b,a+1} :  a \geq b \geq k , a+b =n-1\}$.
This last fact is proved by induction on $n-k$.

It is left to show that the images of $A^S$ are independent in $(\X \oplus \Y) / Z $, 
and span $((\X \oplus \Y)/Z)^{\psi_1,\dots,\psi_6}$ ($\psi_1,\dots,\psi_6$ are indeed well defined on 
$(\X \oplus \Y) /  Z$ since they vanish on $Z$.)
Let $P_Z : \X \oplus \Y \to \Y$ be the projection with respect to the direct summand $Z$.
Then the above is the same as proving that 
$P_Z \circ F^\1(A^S)$ are independent, 
and span $\Y^{\psi_1,\dots,\psi_6}$.

For each $a \in \Z$, $P_Z \circ F^\1(S_{\h{a},0,0})$ is 
$$\gamma_a = -\frac{1}{2}(a-3)Y_{a-3} +\frac{3}{2}(a-1)Y_{a-1} 
 - \frac{3}{2}(a+1)Y_{a+1} + \frac{1}{2}(a+3)Y_{a+3}.$$ 
One checks directly that indeed each $\gamma_a$ satisfies the six equation.
Since these six equations are clearly independent on $\Y$, 
it is enough to show that $\{ \gamma_a \}_{a \in \Z}$ are independent and span a
subspace of $\Y$ of codimension 6. 
For this we show that $\{Y_d \}_{-2 \leq d \leq 2} \cup \{ Y_3 - Y_{-3} \} \cup \{ \gamma_a \}_{a \in \Z}$ are
a basis for $\Y$. 
For this it is enough to prove that for every $n \geq 3$, the $2n+1$ dimensional subspace of $\Y$ 
$span \{ Y_d \}_{-n \leq d \leq n}$, is spanned by the $2n+1$ elements
$\{Y_d \}_{-2 \leq d \leq 2} \cup \{ Y_3 - Y_{-3} \} \cup \{ \gamma_a \}_{-n+3 \leq a \leq n-3}$.
This is shown by induction on $n$.
\end{pf}

From Lemma \pr{gen} and Lemma \pr{inj} we conclude:

\begin{thm}\label{c}
\begin{enumerate}
\item $A^J \cup A^S$ is a basis for $\G$.
\item $F^\1 : \G \to \X \oplus \Y$ is injective, with image $(\X \oplus \Y)^{\psi_1, \dots , \psi_6}$.
\end{enumerate}
\end{thm}

We may now prove our main result:

\begin{thm}\label{main} 
$F(\C_0) \su (\X \oplus \Y)^{\psi_3,\psi_4,\psi_5,\psi_6}$ and 
$F:\C_0 \to (\X \oplus \Y)^{\psi_3,\psi_4,\psi_5,\psi_6}$ is a universal order 1 invariant.
\end{thm}

\begin{pf}
By Proposition \pr{iso} it is enough to show that $\phi_F : \h{\G} \to \X \oplus \Y$ is injective with image
$(\X \oplus \Y)^{\psi_3,\psi_4,\psi_5,\psi_6}$. Indeed, by Theorem \pr{c}, 
$\phi_F|_{\G}$ is injective, and onto $(\X \oplus \Y)^{\psi_1, \dots , \psi_6}$.
It remains to show that $\phi_F(a_0)=F(\Gamma_0)$ and $\phi_F(a_1)=F(\Gamma_1)$ are in 
$(\X \oplus \Y)^{\psi_3,\psi_4,\psi_5,\psi_6}$, and independent mod $(\X \oplus \Y)^{\psi_1, \dots , \psi_6}$.
This follows by observing the values of $\psi_1, \dots , \psi_6$ on $F(\Gamma_0),F(\Gamma_1)$.
We have $F(\Gamma_0) = Y_{-1} + Y_1$ and $F(\Gamma_1) = X_{0,0} + Y_{-2} + Y_0 + Y_2 = Z_{0,0} + Y_{-2} + Y_2$,
so:
\begin{itemize}
\item $\psi_1(F(\Gamma_0))=2$ and $\psi_i(F(\Gamma_0))=0$ for $i \neq 1$. 
\item $\psi_2(F(\Gamma_1))=2$ and $\psi_i(F(\Gamma_1))=0$ for $i \neq 2$. 
\end{itemize}

\end{pf}

From now on we will always present any order 1 invariant in the form $\phi \circ F$.
Whenever this is done in practice, we will define $\phi$ over all $\X \oplus \Y$, which
will allow us to define $\phi$ via the $\{ X_{a,b},Y_d\}$ or the $\{ Z_{a,b},Y_d\}$ basis.
One can then think of $\phi$ as restricted to $(\X \oplus \Y)^{\psi_3,\psi_4,\psi_5,\psi_6}$.

We also remark, that from Theorem \pr{c}(1) it follows that the six families of elements of 
$N$ that we have presented and used 
(Figures \pr{c4},\pr{c5},\pr{c6},\pr{c7} and the mirror images of Figures \pr{c6},\pr{c7}) 
in fact span all $N$. 
Indeed let $N' \su N$ be the subspace of $\F$ spanned by
these elements, and let $K=span(A^J \cup A^S)$. Since the proof of Lemma \pr{gen} used only the relations
coming from these six families, it implies that $K + N' = \F$. 
On the other hand Theorem \pr{c}(1) implies $K \cap N = \{ 0 \}$.
Together this shows $N' = N$.

\section{$S$-invariants and $J$-invariants}\label{siji}

An order 1 invariant $f$ is called an $S$-invariant if $f^\1$ vanishes on all symbols of type $J$.
Similarly, $f$ is called a $J$-invariant if $f^\1$ vanishes on all symbols of type $S$.
We will denote the spaces of (order 1) $S$- and $J$-invariants by $V^S_1$, $V^J_1$ respectively.
Clearly $V^J_1 \cap V^S_1 = V_0$, the space of all invariants which are constant on $\ce$ and on $\co$.
A universal invariant for such subclass of invariants, say $V^S_1$, 
is an $S$-invariant $f^S:\C_0 \to W^S$ such that for any $W$, 
the natural map $Hom_\Q(W^S,W) \to V^S_1(W)$ given by $\phi \mapsto \phi \circ f^S$ 
is an isomorphism. For example, from the proof of Theorem \pr{main} we see that 
$(\psi_1, \psi_2) \circ F : \C_0 \to \Q^2$ is a universal invariant for $V_0$. 
(This is of course nothing other than the invariant assigning $(2,0) , (0,2)$ to all curves in
$\C_0^\od,\C_0^\ev$, respectively.)

We begin with $S$-invariants.
By Theorem \pr{main} any $W$ valued order 1 invariant is of the form $f = \phi \circ F$ for a unique
$\phi : (\X \oplus \Y)^{\psi_3,\psi_4,\psi_5,\psi_6} \to W$. 
By the proof of Lemma \pr{inj}, $f^\1 = (\phi \circ F)^\1 = \phi \circ F^\1$ 
vanishes on all symbols of type $J$  iff 
$\phi$ vanishes on $Z$. So we obtain that $f^S = P_Z \circ F : \C_0 \to \Y^{\psi_3,\psi_4,\psi_5,\psi_6}$ 
is a universal $S$-invariant, where $P_Z$, as in the proof of Lemma \pr{inj}, is the projection
to $\Y$ with respect to the direct summand $Z$. The presentation of $P_Z$ via the $\{ X_{a,b} , Y_d \}$ basis is:
$P_Z(X_{a,b}) = \Phi_{a,b}$, $P_Z(Y_d)=Y_d$.

\begin{example}
Let $\varphi : \X \oplus \Y \to \Q$ be defined via the $\{ Z_{a,b} , Y_d \}$ basis as follows
$\varphi(y_d) = d^2$, $\varphi(Z_{a,b})=0$. By the above analysis $\varphi \circ F$ is an $S$-invariant.
In the $\{ X_{a,b} , Y_d \}$ basis $\varphi$ is given by
$\varphi(y_d) = d^2$, $\varphi(X_{a,b})=4(a-b) - (a+b)^2$. 
By direct substitution (preferably in the $\{ Z_{a,b} , Y_d \}$ basis) we get
$(\varphi \circ F)^\1(S_{\h{a},0,0})  =  24$ for all $a$.
By Proposition \pr{relJS} and Lemma \pr{sk} we conclude that 
$(\varphi \circ F)^\1=  24$ for all symbols of type $S$.
That is, $\frac{1}{24} \varphi \circ F$ coincides with Arnold's Strangeness invariant of spherical curves,
up to choice of constants on $\ce,\co$. As noticed above, such constants may be obtained via
$\psi_1 \circ F, \psi_2 \circ F$, so,  
$(\frac{1}{24} \varphi + k_1 \psi_1 + k_2\psi_2) \circ F$ gives Arnold's strangeness invariant for spherical
curves with all possible normalizations. 
\end{example}

We now look at $J$-invariants. 

\begin{lemma}\label{par}
An order 1 invariant $f$ is a $J$-invariant
iff $f^\1$ vanishes on $A^S$ and
the value of $f^\1$ on the symbols $J^+_{a,b}$ and $J^A_{a,b}$ depends only on the parity of $a$ and $b$. 
\end{lemma}

\begin{pf}
If the 
value of $f^\1$ on $J^+_{a,b}, J^A_{a,b}$ depends only on the parity of $a,b$, then by 
Proposition \pr{relJS} the relations of Lemma \pr{sk} hold, and so if $f^\1$ vanishes on
$A^S$ it vanishes on all symbols of type $S$. The converse is proved similarly,
if $f$ is a $J$-invariant then it vanishes on $A^S$, and by Proposition \pr{relJS}
the value of $f^\1$ on $J^+_{a,b}, J^A_{a,b}$ depends only on the parity of $a$ and $b$ (recall that the indices
$a,b$ are unordered).
\end{pf}

It follows from Lemma \pr{par} that a $J$-invariant  $f : \C_0 \to W$ 
is determined, up to constants, by the value of $f^\1$ on the following six classes of symbols,
where $e,o$ stand for even and odd respectively:
$J^+_{e,e}$, $J^+_{o,o}$, $J^+_{e,o}$, $J^A_{e,e}$, $J^A_{o,o}$, $J^A_{e,o}$.
In addition to the linear functions $\psi_1,\psi_2$ which provide the constants, 
we define the following 6 linear maps $\eta_i : \X \oplus \Y \to \Q$ using the $\{ X_{a,b} , Y_d \}$ basis.
As before, a basis element that is not mentioned is mapped to 0.

\begin{enumerate}

\item $\eta_1(X_{a,b}) = (a+b)^2$ for $a+b$ odd, $\eta_1(y_d) = -d^2$ for $d$ odd.
\item $\eta_2(X_{a,b}) = (a+b)^2$ for $a+b$ even, $\eta_2(y_d) = -d^2$ for $d$ even. 
\item $\eta_3(X_{a,b}) = 1$ for $a$ even $b$ odd.
\item $\eta_4(X_{a,b}) = 1$ for $a$ even $b$ even.
\item $\eta_5(X_{a,b}) = 1$ for $a$ odd $b$ even.
\item $\eta_6(X_{a,b}) = 1$ for $a$ odd $b$ odd.
\end{enumerate}

By direct substitution we see that $\eta_i \circ F^\1$ all satisfy the conditions of Lemma \pr{par}, 
so $\eta_i \circ F$ are all $J$-invariants.
The value of each $\eta_i$ on each of the six classes of symbols of type $J$ mentioned above
is as follows, where blank spaces are 0. 

$$\begin{matrix}
       & J^+_{e,o} &J^A_{e,e} & J^A_{o,o}  & J^+_{e,e} & J^+_{o,o} & J^A_{e,o}   \\  

\eta_1 &    0      &   -8     &   -8             &           &          &                        \\
\eta_3 &      1    &    2     &    0             &           &          &                       \\ 
\eta_5 &      1    &      0   &    2             &           &          &                       \\
 
\eta_2 &           &          &                  &       0   &    0     &   -8              \\
\eta_4 &           &          &                  &       2   &     0    &   1               \\ 
\eta_6 &           &          &                  &       0   &    2     &   1      
\end{matrix}$$

This is a nonsingular matrix, and so, 
if we define $\phi^J : \X \oplus \Y \to \Q^8$ via the eight linear maps
$\psi_1, \psi_2, \eta_1, \dots , \eta_6$, then 
$f^J = \phi^J \circ F : \C_0 \to \Q^8$ is a universal $J$-invariant.

\begin{example}
Let $\phi^+,\phi^-: \X \oplus \Y \to \Q$ be defined by:
\begin{itemize}
\item $\phi^+(X_{a,b}) = 4+(a+b)^2$ for all $a,b$, $\phi^+(Y_d)=-d^2$ for all $d$. 
\item $\phi^-(X_{a,b}) =   (a+b)^2$ for all $a,b$, $\phi^-(Y_d)=-d^2$ for all $d$.
\end{itemize}

Then $\phi^+ = \eta_1 + \eta_2 + 4(\eta_3 + \eta_4 + \eta_5 + \eta_6)$ 
and $\phi^-  = \eta_1 + \eta_2$, so 
$\phi^+ \circ F$ and $\phi^- \circ F$ are $J$-invariants.
By direct substitution or by applying the matrix above, we obtain
that 
\begin{itemize}
\item $(\phi^+ \circ F)^\1(J^+_{a,b}) = 8$,  $(\phi^+ \circ F)^\1(J^A_{a,b}) = 0$, for all $a,b$.  
\item $(\phi^- \circ F)^\1(J^A_{a,b}) =-8$,  $(\phi^- \circ F)^\1(J^+_{a,b}) = 0$, for all $a,b$. 
\end{itemize}

That is, $\frac{1}{4}\phi^+ \circ F, \frac{1}{4}\phi^- \circ F$ 
coincide with Arnold's $J^+,J^-$ invariants of spherical curves,
respectively, up to choice of constants on $\ce,\co$. So, 
$(\frac{1}{4}\phi^+ + k_1 \psi_1 + k_2\psi_2) \circ F$ and $(\frac{1}{4}\phi^- + k_1 \psi_1 + k_2\psi_2)\circ F$ 
give Arnold's $J^+$ and $J^-$ invariants for spherical
curves with all possible normalizations.
\end{example}

An order 1 invariant will be called an $SJ$-invariant if it is the sum of an $S$-invariant and a $J$-invariant.
The space of $SJ$-invariants is thus $V^S_1 + V^J_1$. We
recall that $V^S_1 \cap V^J_1 = V_0$. It follows that 
$f^{SJ} = (P_Z \oplus (\eta_1,\dots,\eta_6)) \circ F : \C_0 \to \Y^{\psi_3,\psi_4,\psi_5,\psi_6} \oplus \Q^6$
is a universal $SJ$-invariant.
We also see that $\phi \circ F$ is an $SJ$-invariant iff $\phi$ vanishes on $Z^{\eta_1,\dots,\eta_6}$.

We have considered several linear maps $H$ on $\X \oplus \Y$ of the special form 
$Y_d \mapsto h(d)$, $X_{a,b} \mapsto -h(a+b)$, for function $h$ on $\Z$,
namely, $\psi_1,\psi_2,\psi_3,\psi_4,\eta_1,\eta_2$ 
and $\phi^-  = \eta_1 + \eta_2$. For $H$ of this form we can interpret $H \circ F(c)$ in terms of the smoothing
of $c$ and its complementary regions in $S^2$. Indeed, from Figure \pr{c1}b it is clear that for $H$ of this form,
$H\circ F (c) = \sum_E \chi(E) h(d(E))$, where the sum is over all complementary regions $E$ of the smoothing of $c$,
and $d(E)$ is minus the sum of winding numbers of the components of the smoothing of $c$, considered in 
$S^2 - \{ p \} \cong \R^2$ for some $p \in E$. (Note that the total winding number is preserved under smoothing.) 
Conversely, this equality shows that for any $h$ the invariant $\sum_E \chi(E) h(d(E))$ is of order 1.
Via this interpretation, there is a clear similarity between our formulas above 
for the two invariants $J^+,J^-$ for spherical curves, and Viro's formulas for $J^+,J^-$ for planar curves, 
appearing in \ct{vi}.
The notion of \emph{index} with respect to a point in $\R^2$, appearing in Viro's formulas, is replaced 
in our formulas by the notion of \emph{winding number} with respect to a point in $S^2$.

\section{Odds, evens, and ends}\label{oee}

In this section we will clarify how all our above results and constructions split between $\ce$ and $\co$.
This should give a better understanding of various details of our work.
We will also obtain two equalities satisfied by any spherical curve, using our linear maps $\psi_5,\psi_6$.

Let $\xe \su \X$ (respectively $\xo \su \X$) 
be the subspace spanned by $X_{a,b}$ with $a+b$ even (respectively odd).
Let $\ye \su \Y$ (respectively $\yo \su \Y$) 
be the subspace spanned by $Y_d$ with $d$ even (respectively odd).

\begin{prop}\label{su1}
$F(\C_0^\ev) \su \xe + \ye$ and
$F(\C_0^\od) \su \xo + \yo$.
\end{prop}

\begin{pf}
Let $c \in \C_0^\ev$. Given a double point $v$ of $c$,
let $U$ be a small neighborhood of $v$, $D= S^2 - U$ and $c_1,c_2$ the two exterior arcs in
$D$. Since the endpoints of $c_1,c_2$
are not intertwined along $\pa D$, there exists a regular homotopy of $c_1,c_2$ in $D$ 
at the end of which they are disjoint. Such regular homotopy preserves $i(c_1),i(c_2)$,
so we may assume $c_1,c_2$ are disjoint. Since $c \in \C_0^\ev$, the number of double points of $c$ is odd
and so the  total number of double points of $c_1$ and $c_2$ is even. Since the parity of $i(c_1),i(c_2)$ coincides
with the parity of the number of double points of $c_1,c_2$, the contribution of $v$ to $f^X(c)$ is in $\xe$.
This is true for each double point $v$, and so $f^X(c) \in \xe$.
The fact that $f^Y(c) \in \ye$ is by definition of $\ce$. The same argument shows that for $c \in \C_0^\od$,
$F(c) \in \xo + \yo$.
\end{pf}

In fact, from the proof of Theorem \pr{main} we obtain the following more precise statement:

\begin{thm}\label{image}
\begin{enumerate}
\item The affine span of $F(\C_0^\ev)$ is the affine subspace of $(\xe + \ye)$
determined by the following three equations on $(\xe + \ye)$:
$\psi_2=2,\psi_4=0,\psi_6=0$.
\item The affine span of $F(\C_0^\od)$ is the affine subspace of $(\xo + \yo)$
determined by the following three equations on $(\xo + \yo)$:
$\psi_1=2,\psi_3=0,\psi_5=0$.
\end{enumerate}
\end{thm}

Let $\F^\ev \su \F$ be the subspace spanned by all symbols
corresponding to curves in $\C_1^\ev$. Let $N^\ev \su N$ be the subspace spanned by 
$v(c)$ for all generic loops in $\ce$, then $N^\ev \su F^\ev$, and let $\G^\ev = \F^\ev / N^\ev$.
Let $\F^\od, N^\od, \G^\od$ be defined in the same way. Then $\F = \F^\ev \oplus \F^\od$
and $\G = \G^\ev \oplus \G^\od$.

\begin{prop}\label{su2}
$F^{(1)}(\G^\ev) \su \xe + \ye$ and 
$F^{(1)}(\G^\od) \su \xo + \yo$. 
\end{prop}

\begin{pf}
For each generating symbol $T$ of $\G^\ev$, $F^{(1)}(T)$ is the difference $F(c)-F(c')$ for two
curves $c,c' \in \C_0^\ev$ (the two resolutions of a curve in $\C_1^\ev$) and so by Proposition \pr{su1},
$F^{(1)}(T) \in \xe + \ye$. The same argument holds for a generating symbol of $\G^\od$.
\end{pf}

Again our work implies a more precise statement:

\begin{prop}\label{su3}
$F^{(1)}(\G^\ev) = (\xe + \ye)^{\psi_2,\psi_4,\psi_6}$ and 
$F^{(1)}(\G^\od) = (\xo + \yo)^{\psi_1,\psi_3,\psi_5}$. 
\end{prop}

Finally, we would like to determine which are the symbols corresponding to curves in $\C_1^\ev$ 
and which to curves in $\C_1^\od$.
We can do this by observing the curves themselves, in which case we must be cautious regarding the way 
the ends of the different exterior curves are intertwined along $\pa D$. Alternatively, we can
use Proposition \pr{su2} and our explicit formulas for $F^\1$.
We obtain: The symbols corresponding to curves in $\C_1^\ev$ are:
$J^+_{a,b}$ for $a+b$ even, $J^A_{a,b},J^B_{a,b}$ for $a+b$ odd, $S_{a,b,c}, S_{\h{a},\h{b},\h{c}}$ for $a+b+c$ even, 
and $S_{\h{a},b,c} , S_{\h{a},\h{b},c}$ for $a+b+c$ odd.
The symbols corresponding to curves in $\C_1^\od$ are the complementary set, that is:
$J^+_{a,b}$ for $a+b$ odd, $J^A_{a,b},J^B_{a,b}$ for $a+b$ even, $S_{a,b,c}, S_{\h{a},\h{b},\h{c}}$ for $a+b+c$ odd, 
and $S_{\h{a},b,c} , S_{\h{a},\h{b},c}$ for $a+b+c$ even.

We conclude this work with the presentation of two equalities satisfied by any spherical curve. 
Define $x_{a,b}:\C_0 \to \Z$ and $y_d : \C_0 \to \Z$ via
$$F(c) = \sum_{a,b} x_{a,b}(c)X_{a,b} + \sum_d y_d(c)Y_d.$$ 
That is, $x_{a,b}(c)$ is the number of double points of $c$ of type $(a,b)$ and $y_d(c)$ is the number of
complementary regions of $c$ of type $d$. 
Theorem \pr{image} gives six equalities satisfied by the curves in $\C_0$,
which can be interpreted as six relations between the invariants $x_{a,b}(c),y_d(c)$.
We exclude the equalities coming from $\psi_1,\psi_2,\psi_3,\psi_4$
since they can be proved directly, by Euler Characteristic arguments, via the observation
in the concluding paragraph of Section \pr{siji}.
The equalities coming from $\psi_5$ and $\psi_6$ are:

\begin{thm}\label{fin}
The following equalities hold for any spherical curve $c \in \C_0$:

\begin{enumerate}

\item $$\sum_{d \ \text{\emph{odd}}} \frac{1}{d}y_d(c) 
+ \sum_{a+b \ \text{\emph{odd}}}  \frac{4(a-b+1) - (a+b)^2}{(a+b)\left( (a+b)^2 - 4 \right)} x_{a,b}(c) = 0.$$

\item $$y_0(c) + \sum_{a+b=0} (b-a-1)x_{a,b}(c) + \sum_{a+b=\pm 2} \frac{a-b}{2} x_{a,b}(c)= 0.$$ 

\end{enumerate}

\end{thm}

Note that by Proposition \pr{su1}, equality (1) is trivially satisfied on $\ce$ and 
equality (2) is trivially satisfied on $\co$.

\end{document}